\documentclass[12pt,a4paper]{amsart}
\usepackage{amsmath}
\usepackage{amssymb}

\usepackage[arrow,matrix]{xy}
    \CompileMatrices 

\usepackage{graphics}
\input xy
\xyoption{all}
\CompileMatrices

\newtheorem{theo}{Theorem}[section]
\newtheorem{remarkk}[theo]{Remark}
\newenvironment{rem}{\begin{remarkk}\rm}{\end{remarkk}}
\newtheorem{definition}[theo]{Definition}
\newenvironment{devery easyfi}{\begin{definition}\rm}{\end{definition}}
\newtheorem{prop}[theo] {Proposition}
\newtheorem{cor}[theo]{Corollary}
\newtheorem{lemma}[theo]{Lemma}
\newtheorem{example}[theo]{Example}

\newtheorem{problem}[theo]{Problem}

\newcommand{\Q}{\ensuremath{\mathbb{Q}}}

\newcommand{\hol}{\ensuremath{\mathcal{O}}}

\newcommand{\Proof}{{\it Proof. }}
\usepackage[arrow,matrix]{xy}
    \CompileMatrices 

\usepackage{graphics}
\input xy
\xyoption{all}
\CompileMatrices

\begin{document}

\title{A volume maximizing canonical surface in $3$-space}
\author{I. C. Bauer}
\author{F. Catanese}
\date{\today}

\maketitle

\section{Introduction}

At the onset of surface theory surfaces in 3-space, and especially
canonical surfaces in 3-space, occupied a central role.

In particular, this study led to the famous Noether inequality $K^2 
\geq 2 p_g - 4$, while
  Castelnuovo observed that  if the canonical map of a minimal smooth 
surface $S$ is birational
(obviously then $p_g \geq 4$) the
inequality $K^2 \geq 3 p_g - 7$ must hold true.

These are the lower bounds for surface geography, but upper bounds 
played a decisive role
in the investigations of the last 30 years, leading to
the socalled Bogomolov-Miyaoka-Yau inequality $$K^2 \leq 9 \chi : = 9 
(p_g -q +1)$$ (cf. \cite{BPV} Chapter
VII, section 4).

For instance, the BMY - inequality gives a hypothetical upper bound for
a question  raised by F. Enriques (cf. \cite{enriques},
chapter VIII, page $284$).

\begin{problem}
Which are the possible values of  $K^2$, in particular which is the 
highest possible value of  $K^2$ for
minimal surfaces
  with geometric genus $p_g = 4$ having a birational canonical map 
(so-called {\em simple canonical}  surfaces)?
\end{problem}

In fact, Enriques even conjectured that the highest possible value for
$K^2$ should be $24$, based on the conjecture that the expected number
of moduli should be strictly positive. The second author showed in 
\cite{bidouble}
that this bound does not hold
true,
constructing simple canonical surfaces with geometric genus $p_g = 4$ 
and $11 \leq K^2 \leq 28$
(in these examples $K^2$ equals the 'canonical degree', i.e., the 
degree of the canonical image).
This bound was improved by C. Liedtke (cf. \cite{liedke}) who showed 
the existence of a simple canonical
surface with $p_g = 4$ and $K^2 = 31$ (and canonical degree $ 12 $).

  Simple canonical surfaces have $K^2
\geq 5$, and for  $5  \leq K^2
\leq 7$ they were constructed by Enriques,
Franchetta, Kodaira, Maxwell, and  for $6 \leq K^2 \leq 16$ by Burniat,
while Ciliberto was able to show for  $5  \leq K^2
\geq 10$ the existence of simple canonical surfaces with ordinary 
singularities (cf. \cite{enriques},
\cite{franchetta},
\cite{maxwell}, \cite{kodaira},
  \cite{burniat}, \cite{ciliberto}).

If we try to go up with $ K^2$,  the BMY inequality tells us that $ 
K^2 \leq 45$, and that,
if equality holds,
then necessarily $S$ is regular $ (q(S) = 0).$

The main result of this paper is the following

\medskip

{\bf Main Theorem.} {\em There exists  a minimal smooth algebraic surface $S$ of
general type over the complex numbers with $K^2 = 45$ and $p_g =
4$, and with birational canonical map. $S$ is rigid, the canonical system
$|K_S|$ has a fixed part and the degree of the canonical image is $19$.}

\medskip

The rigidity of $S$ is due to the fact that, by Yau's proof of the 
inequality $K^2 \leq 9 \chi$,
  it follows (cf. also \cite{miya}) that $K^2 = 9 \chi$
  if and only if the universal covering of $S$ is the 2-dimensional 
complex ball $\mathfrak B_2$.

It was for long time extremely  hard to give direct algebro geometric 
constructions of such ballquotients,
until a breakthrough came via the the explicit constructions by Hirzebruch as
  Kummer coverings of the complex projective plane branched in a
configuration of lines (\cite{hirz}).
  These examples were extended and generalized in
the book \cite{bhh}, which amply describes  three examples of such 
(compact) ballquotients.
The  configurations occurring  are quite classical: a complete 
quadrangle, the Hesse configuration and the dual
Hesse configuration. Even if it is possible to determine the 
numerical data which a configuration has to
fulfill in order to give rise to a ball quotient, it is less easy to compute
the holomorphic invariants.

In fact, already the determination of the irregularities $q$ of the 
Hirzebruch examples
and of some \'etale quotients of them  required
further work  by M.-N. Ishida (cf. \cite {ishida1},
  \cite{ishida2}), but no regular examples were indeed found (except 
Mumford's fake projective  plane,
whose construction however was not so explicit as Hirzebruch's one,
see \cite{mum}).

The example of \cite{bhh} we are interested in here is the 
$(\mathbb{Z}/5 \mathbb{Z})^5$-covering
  $\hat{S}$ of $\mathbb{P}^2$ branched exactly in  a complete 
quadrangle. This surface has the invariants $K^2
= 45 \cdot 125$ and $\chi = 5 \cdot 125$.  It is clear that an 
\'etale $(\mathbb{Z}/5 \mathbb{Z})^3$ quotient
or, equivalently, a smooth $(\mathbb{Z}/5 \mathbb{Z})^2$ covering of 
$\mathbb{P}^2$ branched exactly in a
complete quadrangle has the invariants $K^2 = 45$ and $\chi = 5$. 
Since, as we observed, $\chi = p_g - q +1$,
we have to produce an example of a surface $S$ which is regular 
(i.e., $q=0$) in order to get
the desired  example of a surface with
$K^2 = 45$ and
$p_g =4$. In fact, we will show that up to isomorphisms there are 
exactly four smooth surfaces with $K^2 =
45$, $\chi = 5$, obtained as $(\mathbb{Z}/5 \mathbb{Z})^2$ coverings of 
$\mathbb{P}^2$ branched exactly in a
complete quadrangle: but only one of them is regular (has $q=0$).

\medskip
The main ingredient of our investigation is the theory of Abelian 
Galois coverings,
developped
by Pardini (cf. \cite{Pardini1}), but apparently not sufficiently known.
Since the treatment by Pardini is very algebraic, and at some points
not so explicit, we devote  section $1$ to explain the structure theorem for 
such
Abelian coverings, and especially the relation occurring between the 
topological data (which
allow to construct the examples) and the explicit determination of 
the character sheaves (or eigensheaves) of the covering
(these determine not only  the  topological but also the holomorphic 
invariants of the constructed surface).

Sections  $2$ and $3$ are devoted to the construction of our surfaces, and
to the investigation of the symmetries of our construction.
This study allows us to classify all the examples up to isomorphisms.

\section{Abelian Covers}
In this paragraph we will recall the structure theorem for normal Abelian
Galois ramified coverings. We shall give a more direct presentation 
than the one in the
original paper by R. Pardini (cf. \cite{Pardini1}). This will turn 
out to be more suitable for our purposes.

Let $X$, $Y$ be normal projective varieties, assume $Y$ to be smooth 
and let $\pi : X \rightarrow Y$
be a finite Galois cover with Abelian Galois group $G$. By the theorem
on the purity of the branch locus the critical set of $\pi$ is a 
divisor $R$, the ramification
divisor, whose image $D := \pi (R)$ is called the branch divisor. In 
the case where also
$X$ is smooth we have the following result
(cf. \cite{catmod}, prop. 1.1).

\begin{prop}
If $X$ is smooth $R$ is a normal crossing divisor with smooth components. 
Moreover, if $x \in X$, then
  the stabilizer of $x$ is the direct sum of the stabilizers of the 
components of $R$ passing through $x$ and
these last groups are cyclic.
\end{prop}

We may assume without loss of generality, and will assume in the 
following, that $Y$ is
smooth and $D$ is a normal crossing divisor.

We remark that $\pi$ factors canonically as

$$\xymatrix{
X\ar[rr]^\pi \ar[dr] _p& &Y \\
&X'\ar[ur]_{p'} &
}$$

where $X'$ is maximal such that $p' : X' \rightarrow Y$ is 
unramified. In fact, one takes
  $X' := X / G'$, where $G'$ is the subgroup of $G$ generated by the 
stabilizers $G_x$ of points $x \in X$.

\begin{definition}
$\pi$ is called {\em totally ramified} iff $p'$ is an isomorphism 
(i.e., $G = G'$).
\end{definition}
Observe that $\pi$ is necessarily totally ramified if $Y$ has a 
trivial algebraic fundamental group.

Now, $\pi$ is determined by the surjective homomorphism $\phi : \pi_1 
(Y - D) \rightarrow G$,
  which factors through $\varphi : H_1 (Y - D, \mathbb{Z}) \rightarrow 
G$, since $G$ is assumed to be Abelian.

We denote by $G^*$ the group of characters of $G$, and
we shall use the additive notation for the group
operation in $G^*$ . Recall that $\pi$ 
is flat (for this it
suffices that $Y$ is smooth and $X$ is normal) and that the action of 
$G$ induces a splitting of the direct image of $\hol_X$ into eigensheaves
$$
\pi _* \mathcal{O}_X = \bigoplus_{\chi \in G^*} \mathcal{L}_{\chi}^{-1},
$$
where $G$ acts on the invertible sheaf $\mathcal{L}_{\chi}^{-1}$ via 
the character $\chi$.

  Note that $\mathcal{L}_1 \cong \mathcal{O}_Y$ and denote by 
$L_{\chi}$ a divisor
associated to the eigensheaf $\mathcal{L}_{\chi}$ (thus 
$\mathcal{L}_{\chi}\cong  \hol(L_{\chi})$).

We shall show how:

\begin{itemize}
{\em \item[1)] one calculates $H_1(Y-D, \mathbb{Z})$;

\item[2)] one calculates the character sheaves $\mathcal{L}_{\chi}= 
\hol(L_{\chi})$
   in terms of the surjective homomorphism $\varphi : H_1 (Y - D, 
\mathbb{Z}) \rightarrow G$.}
\end{itemize}

Consider the exact sequence
\begin{equation}\label{exactsequ}
0 \rightarrow K \rightarrow H_1(Y-D, \mathbb{Z}) \rightarrow H_1(Y, 
\mathbb{Z}) \rightarrow 0.
\end{equation}

\begin{rem}\label{tr}
If  $\pi$ is totally ramified, $\varphi|K : K \rightarrow G$ is surjective.
\end{rem}

Following the arguments in \cite{catmod} we obtain

\begin{prop}
$$
K = ker(H_1(Y-D) \rightarrow H_1(Y)) =  coker(r : H^{2n-2}(Y) 
\rightarrow H^{2n-2}(D)).
$$
In particular, if $H_1(Y, \mathbb{Z}) = 0$, then

$H_1(Y-D, \mathbb{Z}) = coker(r : H^{2n-2}(Y) \rightarrow H^{2n-2}(D))$.
  \end{prop}

\Proof
  Let $V$ be an open tubular neighbourhood of $D$ and denote by 
$\partial V$ its boundary.
  Then we have the exact sequence
$$
\ldots \rightarrow H^{2n-2}(Y) \rightarrow H^{2n-2}(D) \rightarrow 
H^{2n-1}(Y,\overline{V})
\rightarrow H^{2n-1}(Y) \rightarrow \ldots
$$

Observing that $H^{2n-1}(Y,\overline{V}) \cong H_1 (Y-D, \mathbb{Z})$, 
we see that
$$
K = ker(H_1(Y-D, \mathbb{Z})\rightarrow H_1(Y, \mathbb{Z})) \cong
$$
$$
\cong coker(r : H^{2n-2}(Y) \rightarrow H^{2n-2}(D) \cong 
\bigoplus_{i=1}^k [D_i] \mathbb{Z}).
$$
\qed

\begin{rem}
Applying $Hom_{\mathbb{Z}}(\_ ~,G)$ to the short exact sequence 
(\ref{exactsequ}) above we get
$$
0 \rightarrow Hom(H_1(Y, \mathbb{Z}),G) \rightarrow Hom(H_1(Y-D, \mathbb{Z}),G)
\rightarrow  Hom(K,G) \rightarrow
$$
$$
\rightarrow Ext^1(H_1(Y, \mathbb{Z}),G) \rightarrow Ext^1(H_1(Y-D, \mathbb{Z}),G)
\rightarrow  Ext^1(K,G) \rightarrow 0.
$$
Therefore an Abelian covering of $Y$ ramified in $D$ is uniquely
determined by a surjective morphism $\varphi : K \rightarrow G$ if and
only if $0 = Hom(H_1(Y),G) $ and $ Ext^1(H_1(Y),G) \rightarrow 
Ext^1(H_1(Y -  D),G)$ is injective. This happens,
for instance, if $H_1(Y, \mathbb{Z})=0$, or more generally if $H_1(Y, 
\mathbb{Z})$ is a finite group whose
exponent
  is relatively prime to the exponent of $G$.
\end{rem}

Let us determine the character sheaves of the Abelian covering determined by
$\varphi : H_1(Y-D, \mathbb{Z}) \rightarrow G$.

Let $\chi \in G^*$ be a character of $G$, i.e., $\chi : G \rightarrow 
C \subset \mathbb{C}^*$,
where $C$ is cyclic. Then $\chi$ induces a surjective morphism $ \chi 
\circ \varphi : H_1(Y-D, \mathbb{Z})
\rightarrow C$, whence a factorization of $\pi$ as

$$\xymatrix{
X\ar[rr]^\pi \ar[dr] & &Y \\
&Z:= X / (ker (\chi \circ \varphi))\ar[ur]_{\pi_{\chi}} &
}$$

where $\pi_{\chi} : Z \rightarrow Y$ is a cyclic covering with group $C$.

\begin{rem}
$\mathcal{L}_{\chi}(Z) = \mathcal{L}_{\chi}(X)$, and we are reduced 
to calculate
the character sheaves for cyclic coverings.
\end{rem}
Write $D = \bigcup_{i=1}^k D_i$ as a union of smooth irreducible components and 
denote by $\delta_i$ the image of a small
loop  around $D_i$  in $H_1(Y-D, \mathbb{Z})$.

Let $d$ be the order of $C$ and let us identify  $C$ with 
$\mathbb{Z} / d$; then we have the well known formula (cf.  the
proof of prop. 4.5 of \cite{torelli})
$$
\hol_Y(d L_{\chi}) \cong \hol_Y(\sum_{i=1}^k (\chi \circ \varphi) 
(\delta_i) \ D_i).
$$

\begin{rem}
We remark that the above linear equivalence
$$
dL_{\chi} \equiv \sum_{i=1}^k (\chi \circ \varphi) (\delta_i) D_i,
$$
  depends only  on $\chi \circ (\varphi|K)$ and does not uniquely determine
the character sheaf $\mathcal{L}_{\chi}$. In fact, if 
$\mathcal{L}_{\chi} \in Pic(Y)$ satisfies the above
equation, then also $\mathcal{L}_{\chi} \otimes \eta$ does, for each 
$d$-torsion sheaf  $\eta \in Pic(Y)$.
If
$\eta$ corresponds to an element $\alpha \in Hom(H_1(Y,\mathbb{Z}), 
\mathbb{Z} / d \mathbb{Z})$, then
$\mathcal{L}_{\chi} \otimes \eta$ is the character sheaf of the 
cyclic covering corresponding to $\chi \circ
\varphi + \alpha \circ p$, where $p : H_1(Y-D, \mathbb{Z}) 
\rightarrow H_1(Y, \mathbb{Z})$. Clearly
$(\chi \circ \varphi + \alpha \circ p)|K = (\chi \circ \varphi)|K$.
\end{rem}

Now, consider the exact sequence

$$
0 \rightarrow \bigoplus_{i=1}^k \mathbb{Z} \cdot D_i \rightarrow 
Pic(Y) \rightarrow Pic(Y-D) \rightarrow 0.
$$

Then $\chi \circ \varphi \in Hom(H_1(Y-D, \mathbb{Z}), \mathbb{Z} / d 
\mathbb{Z})$ corresponds to
  a $d$- torsion sheaf $\eta \in Pic(Y-D)$.  Assume that $\mathcal{L} 
= \mathcal{O}(L)$ is another lifting of
$\eta$ to $Pic(Y)$. Then $\mathcal{L} = \mathcal{L}_{\chi} + \sum a_i 
D_i$, $a_i \in \mathbb{Z}$. Therefore
$$
dL = d(L_{\chi} + \sum a_iD_i) \equiv \sum_{i=1}^k ((\chi \circ 
\varphi) (\delta_i) + da_i) D_i.
$$

Choosing a fixed system of representatives of $\mathbb{Z} / d 
\mathbb{Z}$, e.g.,
  $\mathbb{Z} / d \mathbb{Z} = \{0, \ldots , d-1\}$, we get then the 
uniqueness of $\mathcal{L}_{\chi}$.

\medskip
We will now use  the above approach in order to write explicit equations for $X$ as
a subvariety in the geometric
  vector bundle 
corresponding to the locally free sheaf $\bigoplus_{\chi \in G^*-\{1\}}
\mathcal{L}_{\chi}$.

\begin{rem}

Let $\chi: G \rightarrow C \cong \mathbb{Z} /d$, $\chi ': G 
\rightarrow C \cong \mathbb{Z} /d'$
  be two characters of $G$. Then $ord(\chi + \chi ') = l.c.m.(d,d') 
=:M$. Write $M$ as $M = \lambda \cdot d =
\lambda ' \cdot d'$. Consider the linear equivalences
$$
d L_{\chi} \equiv \sum_{i=1}^k (\chi \circ \varphi) 
(\delta_i) D_i = \sum_{i=1}^k \Delta_i D_i,
$$
$$
d' L_{\chi'} \equiv \sum_{i=1}^k (\chi' \circ \varphi) 
(\delta_i) D_i  = \sum_{i=1}^k \Delta_i' D_i.
$$
and
$$
M ( L_{\chi +\chi'}) \equiv   \sum_{i=1}^k  ((\chi + \chi ')\circ \varphi) 
(\delta_i)  D_i  \equiv   \sum_{i=1}^k ((\lambda \Delta_i + 
\lambda ' \Delta_i ') \ mod (M) )\cdot D_i  \    .
$$

Since  moreover $ 0 < \lambda \Delta_i + \lambda ' \Delta_i ' < 2M$, 
we may write (identifying the divisor $L_{\chi}$ with the divisor
$(\sum_{i=1}^k \frac{\Delta_i}{d} D_i ) \in \oplus_{i=1}^k \Q D_i$)
$$
L_{\chi} + L_{\chi'} - L_{\chi + \chi '} = \sum_{i=1}^k
\epsilon _{D_i}^{\chi, \chi'} D_i,
$$
where $\epsilon _{D_i}^{\chi, \chi'} = 1$ if $\lambda \Delta_i + 
\lambda ' \Delta_i ' \geq M$ and
$\epsilon _{D_i}^{\chi, \chi'} = 0$ otherwise.
\end{rem}

The above equality is equivalent (as shown in  \cite{Pardini1}) to the existence of the 
multiplication maps
$$
\mu_{\chi, \chi '} : \mathcal{L}_{\chi}^{-1} \otimes \mathcal{L}_{\chi'}^{-1}
\rightarrow \mathcal{L}_{\chi + \chi '}^{-1}
$$
which correspond to global sections of $\mathcal{L}_{\chi} \otimes 
\mathcal{L}_{\chi'}
\otimes \mathcal{L}_{\chi + \chi '}^{-1}$ whose divisor
is exactly equal to  $\sum_{i=1}^k
\epsilon _{D_i}^{\chi, \chi'} D_i$.

Let in fact  $\sigma_{i}
\in \Gamma (X, \mathcal{O} (D_{i}))$ a section with $ div ( 
\sigma_{i} ) = D_{i}$: then
$\Pi_{i} \sigma_{i}^{\epsilon_{\chi, \chi 
'}^{i}} $ is a global section of $\mathcal{L}_{\chi} \otimes 
\mathcal{L}_{\chi'}
\otimes \mathcal{L}_{\chi + \chi '}^{-1}$ yielding the multiplication maps.

These sections define  equations for 
the  natural embedding
$$
i : X \hookrightarrow W : = \bigoplus _{\chi \in G^* - \{1\}} 
\mathbb{V} (\mathcal{L}_{\chi}^{-1}).
$$

Let in fact $w_{\chi}$ be a fibre coordinate of $\mathbb{V} 
(\mathcal{L}_{\chi}^{-1})$ :
then $i(X)$ is defined by the equations
$$
w_{\chi} w_{\chi '} = \Pi_{\nu} \sigma_{\nu}^{\epsilon_{\chi, \chi 
'}^{\nu}} w_{\chi+ \chi '}.
$$

We infer the following

\begin{cor}
If $Y$ and the $D_i$'s  are defined over a field $K$, then also $X$ is
defined over  $K$.
\end{cor}

\section{The construction}
We consider in $\mathbb{P}^2 = \mathbb{P}^2_{\mathbb{C}}$ a complete
quadrangle, i.e., the union of the six lines through four points 
$P_0, \cdots P_3$
in general position.

Let $\pi : Y := \hat{\mathbb{P}}^2 (P_0, \cdots , P_3) \rightarrow 
\mathbb{P}^2$
  be the Del Pezzo surface which is the blow up of the plane in 
the points $P_0, \cdots , P_3$.
Denote by
$E_0,
\cdots , E_3$ the exceptional curves. Moreover, for $j = 1, 2, 3$, 
let $L_j' := H - E_0 - E_j$, where $H$ is
the total transform in $Y$ of a line on
$\mathbb{P}^2$, and let $L_j := H - \sum_{i=1}^3 E_i + E_j$. I.e., 
$L_j'$ is the
strict transform of the line in $\mathbb{P}^2$ through $P_0$ and 
$P_j$, whereas $L_j$ is the strict transform
of the line in $\mathbb{P}^2$ through $E_i$ and $E_k$, where 
$\{i,j,k\} = \{1,2,3\}$.

The divisor $L_1 + L_2 + L_3 + L'_1 + L'_2 + L'_3 + E_0 +
E_1 + E_2 + E_3$  on $Y$ has  simple normal crossings and we shall 
denote it by $D$.
\begin{rem}\label{first homology}
It is well known that $H^2(Y, \mathbb{Z})$ is freely generated by 
$H$, $E_0, \cdots , E_3$. Since
\begin{equation}
H_1(Y-D, \mathbb{Z}) \cong coker(r : H^{2n-2}(Y) \rightarrow H^{2n-2}(D) \cong
  \bigoplus_{i=1}^k [D_i] \mathbb{Z}),
\end{equation}
where $r$ is given by the intersection matrix

\bigskip

\begin{tabular}{|l|c c c c c|}
\hline
  & $H$ & $E_0$ & $E_1$ & $E_2$ & $E_3$ \\
\hline
$L_1'$ & 1 & 1 & 1 & 0 & 0\\
$L_2'$ & 1 & 1 & 0 & 1 & 0\\
$L_3'$ & 1 & 1 & 0 & 0 & 1\\
$L_1$ & 1 & 0 & 0 & 1 &  1\\
$L_2$ & 1 & 0 & 1 & 0 &  1\\
$L_3$ & 1 & 0 & 1 & 1 &  0\\
$E_0$ & 0 & -1 & 0 & 0 & 0\\
$E_1$ & 0 & 0 & -1 & 0 & 0\\
$E_2$ & 0 & 0 & 0 & -1 & 0\\
$E_3$ & 0 & 0 & 0 & 0 & -1\\
\hline
\end{tabular}
.

\bigskip

we obtain
$$
H_1(Y-D, \mathbb{Z}) \cong (\bigoplus_{i=0}^3 \mathbb{Z} e_i \oplus 
\bigoplus_{i=1}^3
  \mathbb{Z} l_i \oplus \bigoplus_{i=0}^3 \mathbb{Z} l_i') / H^2(Y, \mathbb{Z}),
$$
where $e_j$ (resp. $l_i$, $l_i'$) is a (small) simple loop oriented 
counterclockwise
around $E_j$ (resp. $L_i$, $L_i'$).
I.e., $H_1(Y-D, \mathbb{Z})$ has generators $e_0, \ldots , e_3, l_1, 
l_2, l_3, l_1',
l_2', l_3'$ and the relations are
$e_0 = l_1' + l_2' + l_3'$, $e_i = l_i' + l_j + l_k$, $\sum l_i' + 
\sum l_i =0$.\\
In particular, $H_1(Y-D, \mathbb{Z})$ is free of rank 5.
\end{rem}

We want to construct a smooth Galois cover $p : S \rightarrow Y$ with 
group $(\mathbb{Z}
  / 5 \mathbb{Z})^2$ branched exactly in $D$.

Such a Galois cover is determined  by a surjective homomorphism

$\varphi : \mathbb{Z}^5 \cong H_1(Y-D, \mathbb{Z} ) 
\rightarrow (\mathbb{Z}
  / 5 \mathbb{Z})^2$ with certain conditions ensuring that $S$ is 
smooth and that the covering branches
exactly in $D$.

By a slight abuse of notation we denote from now on by $\epsilon_h$,
$l_i'$, and $l_j$ the images in $(\mathbb{Z} / 5 \mathbb{Z})^2$  of
the above generators of $H_1(Y-D, \mathbb{Z})$ under the homomorphism 
$\varphi$.

\begin{rem}
It obviously follows from remark (\ref{first homology}) that each 
$\epsilon_h$ is determined by
the $l_i'$, $l_j$ and that $\sum_i l_i' +\sum_j l_j = 0$.
\end{rem}

We write
$$
l_i' :=
\begin{pmatrix}
x_i \\
y_i
\end{pmatrix}
=: u_i; ~~~~~~~~l_j' :=
\begin{pmatrix}
z_j \\
w_j
\end{pmatrix}
=: v_j,
$$
where $x_i$, $y_i$, $z_j$, $w_j \in \{0, \ldots , 4 \} \cong 
\mathbb{Z} / 5 \mathbb{Z}$.

In order to calculate the invariants (i.e., $p_g$, $K_S^2$, $q$) of 
the Galois covering
given by the homomorphism $\varphi$, we have to calculate for each 
character $\chi \in  ((\mathbb{Z} /5
\mathbb{Z})^2)^*$ the eigensheaf  $\mathcal{L}_{\chi}$.

Before doing this let us work out first the two sets of conditions
ensuring that our covering is

1) branched exactly in $D = L_1 + L_2 + L_3 + L'_1 + L'_2 + L'_3 + E_0 +
E_1 + E_2 + E_3$;

and that

2) $S$ is smooth.

\begin{lemma}\label{nec}
1) If $u_i$, $v_j$, $\sum u_i$, $\epsilon_i = u_i + v_j +v_k$ are 
different from zero
in $(\mathbb{Z} /5 \mathbb{Z})^2$, then the covering $p : S 
\rightarrow Y$ is branched exactly in $L_1$,
$L_2$, $L_3$, $L_1$, $L_2'$, $L_3'$, $E_0, \ldots , E_3$.

2) If the following pairs of vectors in $(\mathbb{Z} /5 \mathbb{Z})^2$

$(u_i,v_i)$ for $i \in \{1,2,3\}$, $(u_1,u_1+u_2+u_3)$, $(u_2,u_1+u_2+u_3)$,
  
$(u_3,u_1+u_2+u_3)$, $(u_1, u_1+v_2+v_3)$, $(u_2, u_2+v_1+v_3)$, 
$(u_3, u_3+v_1+v_2)$, $(u_1+ v_2+v_3,v_i)$
for $i=2,3$, $(u_2+ v_1+v_3,v_i)$ for $i=1,3$, $(u_3+ v_1+v_2,v_i)$ for $i=1,2$
are linearly independent, then $S$ is smooth.
\end{lemma}

\Proof 1) is obvious.

2) follows from the fact that $S$ given by the homomorphism $\varphi$ is smooth
if and only if the following condition holds: let $D_1$, $D_2$ be two 
non trivial irreducible subdivisors of
the branch divisor of $p : S \rightarrow Y$ and let $d_i$ a small 
loop around $D_i$, then $\varphi(d_1)$ and
$\varphi(d_2)$ are not in the same cyclic subgroup of $(\mathbb{Z} /5 
\mathbb{Z})^2$. \qed

\begin{rem}
Let $p : S \rightarrow Y$ be a $(\mathbb{Z} /5 \mathbb{Z})^2$ - Galois
cover with $u_i$ and $v_j$ satisfying the two conditions of the above
lemma. Then $S$ is a smooth minimal surface with $K_S^2 = 45$ and
$\chi = 5$. We are interested to find such surfaces with $q=0$,
because then they will have geometric genus equal to $4$.
\end{rem}

Given a character  $\chi = (a,b) \in (\mathbb{Z} /5 \mathbb{Z})^2$,
let us determine $\mathcal{L}_{\chi} = \mathcal{L}_{(a,b)}$.

By the results of section $1$, we get
\begin{prop}\label{character sheaves}
$$
5 \mathcal{L}_{\chi} \equiv \sum_{i=1} ^3 \chi (l_i) L_i +  \sum_{i=1} ^3
\chi (l_i') L_i' + \sum_{i=0} ^3 \chi (e_i) E_i.
$$
$$
5 \mathcal{L}_{(a,b)} \equiv \sum_{i=1} ^3 [ax_i + by_i] L_i + \sum_{i=1} ^3
  [az_i + bw_i] L_i' + [a(x_1+x_2+x_3) +
$$
$$
b(y_1 + y_2 + y_3)] E_0 + \sum_{i=1}^3 [a(x_i + z_j + z_k) + b(y_i + 
w_j + w_k)] E_i.
$$

Here $[z]$ denotes the residue class of $z$ modulo $5$.
\end{prop}

\section{The symmetries of the construction}
\begin{definition}
A six - tuple $\mathfrak{U} := (u_1,u_2,u_3,v_1,v_2,v_3) \in 
((\mathbb{Z} /5 \mathbb{Z})^2 - \{0\})^6$
  is said to be {\em admissible} if and only if

0) $ u_1 + u_2 + u_3 + v_1 + v_2 + v_3 = 0$

and moreover the two conditions of Lemma \ref{nec} are verified:

1) $u_i$, $v_j$, $\sum u_i$, $\epsilon_i = u_i + v_j +v_k$ are 
different from zero in
  $(\mathbb{Z} /5 \mathbb{Z})^2$;

2) the following pairs of vectors in $(\mathbb{Z} /5 \mathbb{Z})^2$

$(u_i,v_i)$ for $i \in \{1,2,3\}$, $(u_1,u_1+u_2+u_3)$, $(u_2,u_1+u_2+u_3)$,
  
$(u_3,u_1+u_2+u_3)$, $(u_1, u_1+v_2+v_3)$, $(u_2, u_2+v_1+v_3)$, 
$(u_3, u_3+v_1+v_2)$, $(u_1+ v_2+v_3,v_i)$
for $i=2,3$, $(u_2+ v_1+v_3,v_i)$ for $i=1,3$, $(u_3+ v_1+v_2,v_i)$ for $i=1,2$,
are linearly independent.
\end{definition}

\begin{rem}
We have seen in the previous section that an admissible six - tuple 
$\mathfrak{U}$
induces a smooth Galois cover $p : S \rightarrow Y$ with Galois group 
$(\mathbb{Z} /5 \mathbb{Z})^2$.
Moreover, $S$ is a minimal surface of general type with $K_S^2 = 45$ 
and $\chi = 5$. We recall that $S$ is a
ball quotient, hence rigid.
\end{rem}

Using MAGMA one sees that there are exactly $201600$ admissible six - tuples.
But of course a lot of them will lead to isomorphic surfaces. In 
order to understand how many non isomorphic
surfaces (with $p_g =4$) we will get by this construction, we have to 
understand the symmetries.

Two admissible six - tuples $\mathfrak{U}$, $\mathfrak{U}'$ obviously 
give isomorphic surfaces
  if there is an automorphism $\phi \in Gl(2, \mathbb{Z}/ 5 
\mathbb{Z})$ such that $\phi (\mathfrak{U}) =
\mathfrak{U}'$. On the other hand the group of biholomorphic 
automorphisms of $\mathbb{P}^2 \setminus
\{L_1,L_2,L_3,L_1',L_2',L_3'\}$ equals $\mathfrak{S}_5$ (cf. 
\cite{Terada}). The action of $\mathfrak{S}_5$ on
the set of admissible six tuples is generated by the following transformations:
$$
(01) : (u_1,u_2,u_3,v_1,v_2,v_3) \rightarrow 
(u_1,u_3+v_1+v_2,u_2+v_1+v_3,u_1+u_2+u_3,v_2,v_3);
$$

$$
(02) : (u_1,u_2,u_3,v_1,v_2,v_3) \rightarrow 
(u_3+v_1+v_2,u_2,u_1+v_2+v_3,v_1,u_1+u_2+u_3,v_3);
$$

$$
(03) : (u_1,u_2,u_3,v_1,v_2,v_3) \rightarrow 
(u_2+v_1+v_3,u_3+v_1+v_2,u_3,v_1,v_2,u_1+u_2+u_3);
$$

$$
(04) : (u_1,u_2,u_3,v_1,v_2,v_3) \rightarrow 
(u_1,u_2,u_3,u_1+v_2+v_3,u_2+v_1+v_3,u_3+v_1+v_2).
$$

It is easy to see that these  four transpositions generate the action 
of a group isomorphic to $\mathfrak{S}_5$.

We consider now the group $\mathcal{G}$ acting on the set of 
admissible six - tuples
  $\mathcal{S}$, which is generated by $\mathfrak{S}_5$ and $Gl(2, 
\mathbb{Z}/ 5 \mathbb{Z})$. Then
$\mathcal{G}$ is a quotient of $Gl(2, \mathbb{Z}/ 5 \mathbb{Z}) 
\times \mathfrak{S}_5$ (the actions commute,
being given by multiplication on the right, respectively on the left).
  A MAGMA computation
shows that $\mathcal{G}$ has four orbits on $\mathcal{S}$. 
Representatives for these orbits are
$$
\mathfrak{U}_1 = (
\begin{pmatrix}
1 \\
0
\end{pmatrix}, \begin{pmatrix}
1 \\
0
\end{pmatrix}, \begin{pmatrix}
0 \\
1
\end{pmatrix}, \begin{pmatrix}
2 \\
1
\end{pmatrix}, \begin{pmatrix}
2 \\
1
\end{pmatrix}, \begin{pmatrix}
4 \\
2
\end{pmatrix});
$$

$$
\mathfrak{U}_2 = (
\begin{pmatrix}
1 \\
0
\end{pmatrix}, \begin{pmatrix}
1 \\
0
\end{pmatrix}, \begin{pmatrix}
0 \\
1
\end{pmatrix}, \begin{pmatrix}
2 \\
1
\end{pmatrix}, \begin{pmatrix}
4 \\
2
\end{pmatrix}, \begin{pmatrix}
2 \\
1
\end{pmatrix});
$$

$$
\mathfrak{U}_3 = (
\begin{pmatrix}
1 \\
0
\end{pmatrix}, \begin{pmatrix}
1 \\
0
\end{pmatrix}, \begin{pmatrix}
0 \\
1
\end{pmatrix}, \begin{pmatrix}
4 \\
1
\end{pmatrix}, \begin{pmatrix}
3 \\
2
\end{pmatrix}, \begin{pmatrix}
1 \\
1
\end{pmatrix});
$$

$$
\mathfrak{U}_4 = (
\begin{pmatrix}
1 \\
0
\end{pmatrix}, \begin{pmatrix}
1 \\
0
\end{pmatrix}, \begin{pmatrix}
0 \\
1
\end{pmatrix}, \begin{pmatrix}
1 \\
1
\end{pmatrix}, \begin{pmatrix}
0 \\
3
\end{pmatrix}, \begin{pmatrix}
2 \\
0
\end{pmatrix}).
$$

The orbit of $\mathfrak{U}_1$ has length $28800$, whereas the orbits 
of $\mathfrak{U_2}$,
$\mathfrak{U_3}$, $\mathfrak{U_4}$ have respective length $57600$.

In particular we see that
$\mathcal{G} \cong
Gl(2, \mathbb{Z}/ 5 \mathbb{Z}) \times \mathfrak{S}_5$.

\medskip
We have moreover:

\begin{theo}\label{four}
Let $S_i$ be the minimal smooth surface of general type with $K^2 
=45$ and $\chi =5$ obtained
  from the covering induced by the admissible six - tuple 
$\mathfrak{U}_i$, where $i \in \{1,2,3,4\}$. Then we
have that $S_3$ is regular (i.e., $q(S_3) = 0$), whereas $q(S_i) = 2$ 
for $i \neq 3$.
\end{theo}

In particular, $S_3$ is the unique minimal surface with $K_S^2 = 45$
and $p_g = 4$ obtained as a $(\mathbb{Z}/ 5 \mathbb{Z})^2$ - cover of
$\mathbb{P}^2$ branched exactly in a complete quadrangle of the
complex projective plane. \\

\Proof
We will calculate the geometric genus of $S = S_3$, using the formula
$$
H^0(S,\mathcal{O}_S (K_S)) = \bigoplus_{(a,b) \in G^*} H^0(Y, 
\mathcal{O}_Y(K_Y) \otimes \mathcal{L}_{(a,b)}).
$$
Applying proposition (\ref{character sheaves}) we obtain the 
following table for
  the character sheaves $\mathcal{L}_{(a,b)}$: \\
{\Small
\begin{tabular}{|l|c|c|c|c|c|}
\hline
$L_{(a,b)}$ & $a = 0$ & $a = 1$  & $a = 2$  & $a = 3$ & $a = 4$ \\
\hline
$b = 0$ & $\mathcal{O}_Y$ & $2H - E_1 - E_2 $ & $2H - E_1 - E_2$ & 
$3H - E_0 - 2E_1$ & $3H - E_0 -2E_1$\\
& &  $ - E_3$ & & $- E_2 - E_3$ &  $- E_2$ \\
\hline
$b=1$ & $H$ & $H$ & $3H - E_0 -E_1$ & $3H -E_0 -E_1$ & $3H -E_0 -E_1$ \\
& & & $-E_2 - E_3$ & $-2E_2 - E_3$ & $-E_2 - E_3$ \\
\hline
$b=2$ & $2H - E_1 - E_3$ & $2H - E_1 - E_2$ & $2H - E_0 - E_1$ & $3H 
-E_0 -E_1$ & $3H - 2E_0 -E_1$ \\
& & $-E_3$ & $-E_2$ & $-E_2 - E_3$ & $-E_2 - E_3$ \\
\hline
$b=3$ & $2H - E_2 - E_3$ & $3H -E_0 -E_1$ & $2H - E_0 - E_3$ & $2H - 
E_0$ & $4H - 2E_0 - E_1$ \\
& & $-E_2 - E_3$ & & & $- 2E_2 - 2E_3$ \\
\hline
$b=4$ & $3H - E_1 - E_2$ & $2H - E_0 -E_3$ & $3H -E_0 -E_1$ & $3H 
-2E_0 -E_1$ & $3H -2E_0 -E_1$ \\
& $-2E_3$ & & $-E_2 - 2E_3$ & $-E_2 - E_3$ & $-E_2$ \\
\hline
\end{tabular}
}

We see immediately that $H^0(Y, \mathcal{O}_Y(K_Y) \otimes
\mathcal{L}_{(a,b)}) = 0$ for all $(a,b) \notin \{(2,1),(3,2),(1,3),
(4,1) \}$ and $H^0(Y, \mathcal{O}_Y(K_Y) \otimes \mathcal{L}_{(a,b)})
\cong \mathbb{C}$ for $(a,b) \in \{(2,1),(3,2),(1,3), (4,1) \}$, i.e.,
$p_g(S_3) = 4$. This proves the claim for $S_3$.

The geometric genus of the remaining surfaces is calculated in
exactly the same way. \qed

\section{The canonical map}
In the previous section we have constructed a minimal surface $S$ of general
type with $K_S^2 = 45$, $p_g = 4$ and $q(S) =0$. We want
now to understand the behaviour of the
canonical map of $S$.

For $(a,b) \in (\mathbb{Z} / 5)^2$ we write
$$
\delta_i(a,b) := [ax_i + by_i]
$$
$$
\lambda_j(a,b) := [az_j + bw_j]
$$
$$
\mu_0(a,b) := [a(x_1+x_2+x_3) + b(y_1 + y_2 + y_3)],
$$
$$
\mu_h(a,b) := [a(x_h+z_j+z_k) + b(y_h + w_j + w_k)].
$$
Then we know that
$$
5 \mathcal{L}_{(a,b)} = \sum_{i=1}^3 \delta_i (a,b)L_i' + 
\sum_{j=1}^3 \lambda_j (a,b)) L_j
  + \sum_{h=0}^3 \mu_h (a,b)) E_h.
$$

Denote by $R_1, \ldots , R_{10}$ the ramification divisors of $p : S 
\rightarrow Y$
lying over $L_1'$, $L_2'$, $L_3'$, $L_1$, $L_2$, $L_3$, $E_0, \ldots 
, E_3$: it is easy to see that they are all irreducible genus 2 curves. Moreover,
  let $x_i$ be a local
equation of $R_i$. We already saw that  $H^0(S, \mathcal{O}_S (K_S))$
is the direct sum of 4 one dimensional eigenspaces $H^0(S, 
\mathcal{O}_S (K_S))_{(a,b)} \cong H^0((Y,
\mathcal{O}_Y(K_Y) \otimes \mathcal{L}_{(a,b)} ) \cong  H^0 ( \mathbb{P}^2, 
\hol_{\mathbb{P}^2})$.
Then a basis of $H^0(S, \mathcal{O}_S (K_S))$ is given by
$$
\{x_1^{4 - \delta_1 (a,b)}\cdot  x_2^{4 - \delta_2 (a,b)} \cdot 
x_3^{4 - \delta_3 (a,b)} \cdot
x_4^{4 - \lambda_1 (a,b)} \cdot x_5^{4 - \lambda_2 (a,b)}\cdot x_6^{4 
- \lambda_3 (a,b)}\cdot
$$
$$
\cdot x_7^{4 - \mu_0 (a,b)} \cdot \cdots x_{10}^{4 - \mu_3 (a,b)} | H^0((Y,
\mathcal{O}_Y(K_Y) \otimes \mathcal{L}_{(a,b)} \neq 0 \}.
$$

It is easy to compute the  table giving the numbers  $\lambda_i$,
$\delta_j$ and $\mu_h$ for $(a,b) \in \{(2,1),(3,2),(1,3), (4,1) \}$: \\

\begin{tabular}{|l|c|c|c|c|c|c|c|c|c|c|}
\hline
$(a,b)$ & $\delta_1$ & $\delta_2$ & $\delta_3$ & $\lambda_1$ &
$\lambda_2$ & $\lambda_3$ & $\mu_0$ & $\mu_1$ & $\mu_2$ & $\mu_3$ \\
\hline
$(1,3)$ & 1 & 1 & 3 & 2 &  4 & 4 & 0 & 4 & 2 & 4 \\
\hline
$(2,1)$ & 2 & 2 & 1 & 4 &  3 & 3 & 0 & 3 & 4 & 3 \\
\hline
$(3,2)$ & 3 & 3 & 2 & 4 &  3 & 0 & 3 & 1 & 2 & 4 \\
\hline
$(4,1)$ & 4 & 4 & 1 & 2 & 4 & 0 & 4 & 3 & 1 & 2 \\
\hline
\end{tabular}

\medskip
Therefore we have the following result

\begin{lemma}\label{basis}
A basis for $H^0(S, \mathcal{O}_S(K_S))$ is given by
$$
\{x_1^3 x_2^3 x_3 x_4^2 x_7^4 x_9^2, x_1^2 x_2^2 x_3^3 x_5 x_6 x_7^4 
x_8 x_{10},
x_1 x_2 x_3^2 x_5 x_6^4 x_7 x_8^3 x_9^2, x_3^3 x_4^2 x_6^4 x_8 x_9^3 x_{10}^2\}.
$$
\end{lemma}

We can now  prove  the following
\begin{theo}
1) The canonical map $\phi_K$ of $S$ has $R_3$ as fixed part  and
its movable part has five base points.
  We have $R_3^2 = -1$ and $K_S.R_3 = 3$. The base points of $K_S - 
R_3$ are $x_1 \cap x_4$ (of
type $(1,1))$, $x_1 \cap x_8$ (of type $(1,1,1)$), $x_2 \cap x_9$ (of 
type $(2,1,1)$), $x_3 \cap x_7$ (of type
$(2,1,1)$), $x_6 \cap x_9$ (of type $(1,1)$).

2) The canonical map is birational and its image in $\mathbb{P}^3$ 
has degree $19$.
\end{theo}

In order to keep the formulation of the above theorem as simple as 
possible we adopted
  the notation: a base point $p$ of $K-R_3$ on $S$ is {\em of type} 
$(n_1, n_2, \ldots , n_k)$ iff $p$ is a
$n_1$ - tuple base point of $K-R_3$, after one blow up the strict 
transform of $K - R_3$ has a $n_2$ - tuple base
point and so on.

\Proof
1) It is immediate from the description of the basis of $H^0(S, 
\mathcal{O}_S(K_S))$
given in lemma (\ref{basis})
that $R_3$ is the only fixed part of $|K_S|$and that $|K-R_3|$ does 
not have a fixed part. It is easy to see
that the base points of $|K-R_3|$ are exactly $x_1 \cap x_4$, $x_1 
\cap x_8$, $x_2 \cap x_9$, $x_3 \cap x_7$,
$x_6 \cap x_9$. Next we will see which kind of base points we have 
and whether there are still infinitely near
base points.

1) $x_1 \cap x_4$: locally around this point $|K-R_3|$ is given by 
$x_1^3x_4^2, x_1^2, x_1, x_4^2$.
  The ideal generated is the ideal $(x_1 , x_4^2)$, thus $K-R_3$ has a 
base point of
type $(1,1)$ in $x_1 \cap x_4$.

2) $x_1 \cap x_8$: locally around this point $|K-R_3|$ is given by 
$x_1^3, x_1^2 x_8, x_1 x_8^3, x_8$.
  The ideal generated is the ideal $(x_1^3 , x_8)$, thus
  $K-R_3$ has a base point of type $(1,1,1)$ in $x_1 \cap x_8$.

Similarly  in the three remaining cases we see that: \\
3) $|K-R_3|$ has a base point of type $(2,1,1)$ in $x_2 \cap x_9$ 
(ideal $(x_2^2 , x_2 x_9^2)$).

4) $|K-R_3|$ has a base point of type $(2,1,1)$ in $x_3 \cap 
x_7$ (ideal $(x_3^2 ,  x_7^4, x_3 x_7)$).

5) $|K-R_3|$ has a base point of type $(1,1)$ in $x_6 \cap x_9$ 
(ideal $(x_6 ,  x_9^2)$).

Therefore we get $deg \phi_K deg \phi_K(S) = (K-R_3)^2 - 2 \cdot 4 - 11 = 19$.
Here we use that $R_3$ has self
intersection $-1$ and genus $2$. It follows immediately that $deg \phi_K = 1$
  and that the canonical image has
degree $19$. \qed

The following

\begin{cor}
There exist surfaces of general type $S$ with birational canonical 
map such that
the canonical system of each deformation of $S$ has base points.
\end{cor}

is an answer to a question by Junho Lee.

\medskip

{\bf Acknowledgement.} We are grateful to Fritz Grunewald for his 
invaluable help
in the explicit calculations leading to Theorem \ref{four}.

This research was performed in the realm of the DFG Schwerpunkt
 "Globale methoden in der komplexen Geometrie".

\vfill

\noindent
{\bf Author's address:}

\bigskip

\noindent
Prof. Dr. Ingrid Bauer \\
  Mathematisches Institut\\
Universit\"at Bayreuth, NWII\\
   D-95440 Bayreuth, Germany

e-mail: Ingrid.Bauer@uni-bayreuth.de

\bigskip
\noindent
Prof. Dr. Fabrizio Catanese\\
Lehrstuhl Mathematik VIII\\
Universit\"at Bayreuth, NWII\\
   D-95440 Bayreuth, Germany

e-mail: Fabrizio.Catanese@uni-bayreuth.de

\end{document}